\documentclass{amsart}
\usepackage[a4paper,margin=32mm]{geometry}
\usepackage{math}
\usepackage[ruled,vlined]{algorithm2e}
\usepackage{mathabx,mathrsfs,mathtools,amsrefs}

\newtheorem{mainthm}{Theorem}
\newtheorem{maincor}[mainthm]{Corollary}

\begin{document}
\title{On commutator length in free groups}
\author{Laurent Bartholdi}
\address{Mathematisches Institut, Universit\"at des Saarlandes}
\email{laurent.bartholdi@gmail.com}

\author{Danil Fialkovski}
\email{19fdr97@gmail.com}

\author{Sergei O. Ivanov} 
\address{
Laboratory of Modern Algebra and Applications,  St. Petersburg State University, 14th Line, 29b,
Saint Petersburg, 199178 Russia}
\email{ivanov.s.o.1986@gmail.com}

\thanks{This work is supported by the Ministry of Science and Higher Education ofor the Russian Federation, agreement 075-15-2019-1619.}

\date{October 11th, 2021}

\begin{abstract}
  Let $F$ be a free group. We present for arbitrary $g\in\N$ a \textsc{LogSpace} (and thus polynomial time) algorithm that determines whether a given $w\in F$ is a product of at most $g$ commutators; and more generally an algorithm that determines, given $w\in F$, the minimal $g$ such that $w$ may be written as a product of $g$ commutators (and returns $\infty$ if no such $g$ exists). The algorithm also returns words $x_1,y_1,\dots,x_g,y_g$ such that $w=[x_1,y_1]\cdots[x_g,y_g]$.

  The algorithms we present are also efficient in practice. Using them, we produce the first example of a word in the free group whose commutator length \emph{decreases} under taking a square. This disproves in a very strong sense a conjecture by Bardakov.
\end{abstract}
\maketitle

\section{Introduction}
Let $F$ be a free group and $[F,F]$ its derived subgroup; so every element $w\in [F,F]$ is a product of commutators $[u,v]=u^{-1}v^{-1}u v$. The minimal number of terms in such a product is called the \emph{commutator length} of $w$. This ``norm'' $\|\cdot\|$ on $[F,F]$ was the subject of much investigation, already by Burnside~\cite{burnside:theory}*{\S238, Ex.~7}, but is still poorly understood, in particular in relation to the usual word length $|w|$. One surprising phenomenon is that $\|w^m\|$ can be smaller than $m\cdot\|w\|$; for $F=\langle x,y\rangle$, we have
\[\|[x,y]^3\|\le2\text{ since }[x,y]^3=[x^{-1}yx,x^{-2}yxy^{-1}]\cdot[yxy^{-1},y^2]\]
(and in fact $=2$; more generally, $\|[x,y]^m\|=\lfloor\frac m2\rfloor+1$, see~\cite{culler:eqs}*{Example~2.6}. In contrast, \emph{stable commutator length}, the limit $\operatorname{scl}(w)=\lim_{m\to\infty}\|w^m\|/m$, is much better understood, see~\cite{calegari:scl}).

\subsection{Algorithms}
The first algorithm for computing commutator length was constructed by Goldstein and Turner~\cite{goldstein-turner:graphtheory}; see also~\cites{culler:eqs,olshanskii:homo}. The method is fundamentally topological: construct a graph with $w$ labeled along a cycle, and determine the minimal genus of a topological surface on which the graph embeds. Given a graph, it is straightforward to compute the minimal genus by linear algebra, but a large number of graphs need to be considered.

Bardakov suggests a more algebraic algorithm in~\cite{bardakov:cl}, see~\S\ref{ss:bardakov}, which translates the problem into a calculation in the symmetric group $S_{|w|}$. He proves that $\|w^m\|$ increases at least linearly with $m$, giving a lower bound based on a quasi-homomorphism, and conjectures $\|w^m\|\ge\frac{m+1}{2}\|w\|$ for all $w,m$. Note that from~\cite{schutzenberger:nmp} we have $\|w^m\|\ge2$ for all $m\ge2$ and $w\neq1$.

Yet a different algorithm is proposed by Calegari~\cite{calegari:scl}*{p.~96sqq}, based on linear programming. In fact, \emph{stable} commutator length can be computed as the solution of a linear program of polynomial size in $|w|$, and integer solutions lead to commutator length. Integer linear programming is much more computationally intensive than real linear programming, but is nevertheless feasible, and has been implemented by Walker, see~\texttt{scallop}~\cite{walker:scallop}.

Our main result is:
\begin{mainthm}\label{thm:main}
  Consider a non-trivial word $w\in [F,F]$ in a free group $F=\langle S\rangle$. Then there exists a factorization of $w$ without cancellations
  \[w=w_1a^{-1}w_2b^{-1}w_3a w_4b w_5\]
  with $a,b\in S\cup S^{-1}$ and $\|w_1w_4w_3w_2w_5\|=\|w\|-1$. Furthermore we have
  \[w=[(w_4 w_3 a)^{w_1^{-1}},(b w_2^{-1} w_3^{-1})^{w_4^{-1}w_1^{-1}}]\cdot(w_1w_4w_3w_2w_5).\]
\end{mainthm}

Recall that a decision problem is in \textsc{LogSpace} if it can be solved by a Turing machine with read-only input and one auxiliary read-write tape initially empty, with the guarantee that its read-write head remains within $\mathcal O(\log n)$ steps of the origin, on an input word of length $n$. Its number of total configurations is bounded by a polynomial in $n$, so such a machine stops after polynomial time if it ever stops. Lipton and Zalcstein prove in~\cite{lipton-zalcstein:wplog} that the word problem in free groups is in \textsc{LogSpace}. From Theorem~\ref{thm:main} we deduce:
\begin{maincor}\label{cor:logspace}
  Let a free group $F$ and an integer $g\in\N$ be fixed. Then the problem ``given $w\in F$, is $\|w\|\le g$?'' is in \textsc{LogSpace}.
\end{maincor}

Let $F$ be a free group. We call \emph{$F$-RAM machine} the extension of the computational model of RAM machines with finitely many registers holding elements of $F$, which can be left- and right-multiplied by generators and tested on their left-most and right-most letter in constant time.
\begin{maincor}\label{cor:algo}
  Let a free group $F$ be fixed. Then there is an algorithm for an $F$-RAM machine that, given $w\in F$, determines $\|w\|$ in time $\mathcal O(|w|^{4\|w\|})$.

  Furthermore, this algorithm returns a representation of $w$ as a product of $\|w\|$ commutators.
\end{maincor}

An algorithm by Wicks~\cite{wicks:commutators} determines whether a word $w$ is a commutator by the following criterion: ``some cyclic permutation of $w$ must be of the form $w_1^{-1}w_2^{-1}w_3^{-1}w_1w_2w_3$ as a product without cancellation''. This leads to a an $\mathcal O(|w|^3)$-time algorithm by searching for the possible starting positions of $w_1,w_2,w_3$. We give, in~\S\ref{ss:implementation}, an algorithm that, assuming constant-time arithmetic operations on integers, improves the average time complexity to $\mathcal O(|w|^2)$.

Given a solution $w=[u_1,v_1]\cdots[u_g,v_g]$ to the problem of expressing $w$ as a product of $g$ commutators, numerous other solutions may be derived by elementary transformations, such as ``replace $u_i$ by $v_i u_i$''. Viewing a solution as a homomorphism $\Sigma_{g,1}\to F$ with $\Sigma_{g,1}=\langle u_1,v_1,\dots,u_g,v_g,c\mid[u_1,v_1]\cdots[u_g,v_g]=c\rangle$ the fundamental group of a surface of genus $g$ with one boundary component, we see that the mapping class group of that surface (which coincides with the outer automorphism group of $\Sigma_{g,1})$ naturally acts on the space of solutions, by precomposition. Culler proves in~\cite{culler:eqs}*{Theorem~4.1} that there are finitely many orbits of solutions under the mapping class group action. We believe that our algorithm produces at least one solution in every orbit.

For example, running our algorithm on $w=x^{-1}y^{-1}x^2yx^{-1}$ with $g=1$ produces, with the convention $x^{-1}=X$ and $y^{-1}=Y$, the solutions
\[w = [x, yX] = [xx, yX] = [Yx, YX] = [yx, XX] = [yx, X]\]
while running it on $w=x^{-1}x^{-1}y^{-1}x^{-1}yxy^{-1}x^2y$ with $g=2$ produces the solutions
\begin{alignat*}{2}
  w &= [YYXyx, YxYxyy]\cdot [y, xy]&
  &= [YYXyx, YxYxyy]\cdot [y, Yxy]\\
  &= [YxYXyx, YxxYxyXy]\cdot [Xy, y]&
  &= [YxYXyx, YxxYxyXy]\cdot [YXy, y]\\
  &= [Yxyx, YxxyxYXy]\cdot [y, X]&
  &= [x, Yxyx]\cdot [x, y]\\
  &= [YXyxYx, xyXYxy]\cdot [y, x]&
  &= [yxYx, xyxyXY]\cdot [Y, X]\\
  &= [xx, xyXYxy]\cdot [y, x]&
  &= [XyxYxx, yyXYx]\cdot [x, Y]\\
  &= [xx, yXYxy]\cdot [y, x]&
  &= [XYxyxx, XYxyx]\cdot [x, y]\\
  &= [x, XYxyx]\cdot [x, y]&
  &= [XXyxYxxx, XyXyXYxx]\cdot [x, Yx]\\
  &= [XXyxYxxx, XyXyXYxx]\cdot [x, XYx]&
  &= [xx, XyXYxy]\cdot [y, x]\\
  &= [XXyxx, x]\cdot [Yx, xy]&
  &= [XXyxx, x]\cdot [xx, xy]\\
  &= [XXyxx, x]\cdot [xx, y]&
  &= [XXyxx, x]\cdot [xx, Xy]\\
  &= [XXYxyxx, XXYxx]\cdot [Yx, xy]&
  &= [XXYxyxx, XXYxx]\cdot [xx, xy]\\
  &= [XXYxyxx, XXYxx]\cdot [xx, y]&
  &= [XXYxyxx, XXYxx]\cdot [xx, Xy]\\
  &= [XXYxxyxx, XXYxyXYXyxx]\cdot [x, yx]&
  &= [XXYxxyxx, XXYxyXYXyxx]\cdot [x, Xyx]\\
  &= [XXYxxyxx, XXXYXyxx]\cdot [x, yx]&
  &= [XXYxxyxx, XXXYXyxx]\cdot [x, Xyx].
\end{alignat*}

\subsection{Non-monotonicity}
We already noted that $\|w^m\|<m\cdot\|w\|$ can occur, for example with $w=[x,y]$. Bardakov's conjecture `$\|w^m\|\ge\frac{m+1}2\|w\|$' is refuted from general results on solutions of equations in free groups: Kharlampovich and Myasnikov deduce~\cite{kharlampovich-myasnikov:genus}*{Theorem~3} the existence of a sequence $(w_n)$ in $[F,F]$ with $\|w_n\|\to\infty$ and $\|w_n^2\|$ bounded; equivalently, the infinite product $F^\infty$ contains an element of order $2$ in its abelianization.

These results are however fundamentally not constructive; and explicit elements with $\|w^2\|<\|w\|$ had eluded discovery (see~\cite{comerford-lee:2comm} and unpublished experiments by Spellman using \texttt{Magma}) up to now. We shall see:
\begin{mainthm}\label{thm:32}
  There exists an element $w\in [F,F]$ of length $64$, see~\eqref{eq:w}, with $\|w\|=3$ and $\|w^2\|=2$.
\end{mainthm}

We obtain this result by running our algorithm on an enumeration of the solutions to a quadratic equation, see~\S\ref{ss:quadratic}. Some tricks and techniques were necessary to render such a computation (and discovery) feasible, which we record in~\S\ref{ss:implementation}.

\section{A commutator length algorithm}\label{ss:commlength}
We begin by recalling standard notation. Fix a free group $F=\langle x_1,\dots,x_r\rangle$. Its elements are represented as \emph{words} over the \emph{letters} $x_1,\dots,x_r,x_1^{-1},\dots,x_r^{-1}$; words may be simplified by removing a pair of contiguous letters $a a^{-1}$ as often as possible, and two words are deemed equal if they simplify to the same word, called \emph{freely reduced}. A word is called \emph{cyclically reduced} if furthermore its first and last letters do not cancel. We denote by $w[i]$ the $i$th letter of a word $w$, numbered from $1$ to $|w|$; and for $j\ge i-1$ we denote by $w[i:j]$ the possibly-empty subword $w[i]w[i+1]\cdots w[j]$. The following is our main algorithm computing commutator length; some improvements will be given in~\S\ref{ss:implementation}:

\LinesNumbered
\begin{algorithm}[H]
  cyclically reduce $w$\;
  \If{$g=0$}{\Return $w=1$}
  \For{$1\le i<j<k<\ell\le|w|$}{
    \If{$w[i]\ne w[k]^{-1}$ \textbf{\textup{or}} $w[j]\ne w[\ell]^{-1}$}{\textbf{continue}}
    \If{$\|w[1:i-1]w[k+1:\ell-1]w[j+1:k-1]w[i+1:j-1]w[\ell+1:|w|]\|\le g-1$}{\Return \textsf{true}}
  }
  \Return \textsf{false}
  \caption{Test recursively whether $\|w\|\le g$}\label{algo:1}
\end{algorithm}

\subsection{Bardakov's theorem}\label{ss:bardakov}
Let us denote by $S_n$ the symmetric group on $n$ elements. We compose elements of $S_n$ right-to-left, as functions. Every $\sigma\in S_n$ induces a partition of $\{1,\dots,n\}$ into orbits. Let $(1,\dots,n)$ denote the cyclic permutation of $\{1,\dots,n\}$, and for a permutation $\pi\in S_n$ write
\[v(\pi)\coloneqq\text{number of orbits of }(1,\dots,n)\pi.\]

Let now $w=w_1\dots w_n$ be a word (not necessarily reduced) over an alphabet $\{x_1^{\pm 1}, \dots, x_r^{\pm 1} \}$. A \emph{pairing} on the word  $w$ is an involution $\pi \in S_n$ such that  $a_{\pi(i)}=a_i^{-1}$ for all $i$. Note that a pairing exists if and only if $w$ represents an element of the commutator subgroup.
\begin{thm}[\cite{bardakov:cl}*{Theorem~1}]\label{thm:bardakov}
  Consider a word $w$ representing a non-trivial element also written $w$ of the derived subgroup. Then
  \[\|w\|=\min_{\pi\text{ pairing on }w}\left(\frac{1-v(\pi)}2+\frac n4\right).\]
\end{thm}
\begin{rem}
  Bardakov proves his result only in the case of cyclically reduced words, but does not use this in the proof. He also overlooks the restriction that $w$ be non-trivial.
\end{rem}

\subsection{Induced permutations}
Let us consider more generally the symmetric group $S_X$ on a set $X$, and an injective map $\alpha\colon Y\to X$. There is a \emph{renormalization map} $S_X\to S_Y$ induced by $\alpha$, and defined as follows: for $\sigma\in S_X$ and $y\in Y$ let $m(y)\ge1$ be the least positive integer such that $\sigma^{m(y)}(\alpha(y))\in\alpha(Y)$, and set
\[\sigma_\alpha(y)\coloneqq \alpha^{-1}(\sigma^{m(y)}(\alpha(y))).\]
Note that the renormalization map is not quite a homomorphism; nevertheless, we have
\begin{lem}\label{lem:if1}
  Let $\alpha\colon Y\to X$ be an injective map, and let $\tau\in S_X$ be such that $\tau(\alpha(Y))=\alpha(Y)$, so $\tau=\tau'\tau''$ with $\tau'$ fixing $\alpha(Y)$ and $\tau''$ fixing $X\setminus\alpha(Y)$. Then for all $\sigma\in S_X$ we have
  \[(\sigma\tau)_\alpha=(\sigma\tau')_\alpha\tau_\alpha.\]
\end{lem}
\begin{proof}
  We first note $\tau''_\alpha=\alpha^{-1}\tau\alpha=\tau_\alpha$. Consider then $y\in Y$, write $x\coloneqq\alpha(y)$, and let $m>0$ be minimal such that $(\sigma\tau)^m(x)\in\alpha(Y)$; thus $(\sigma\tau)_\alpha(y)=\alpha((\sigma\tau)^m(x))$. Now $\tau(x)=\tau''(x)$ by assumption, while $z_i\coloneqq(\sigma\tau)^i(x)\notin\alpha(Y)$ for $i<m$ so $\tau(z_i)=\tau'(z_i)$; thus $(\sigma\tau)^i(x)=(\sigma\tau')^i(\tau''(x))$, and $m$ is also minimal such that $(\sigma\tau')^m(\tau''(x))\in\alpha(Y)$. 
\end{proof}

\begin{lem}\label{lem:if2}
  Let $\alpha\colon Y\to X$ be an injective map and consider $\sigma\in S_X$. Then the map $O\mapsto \alpha^{-1}(O)$ is a bijection between the orbits of $\sigma$ on $X$ that intersect $\alpha(Y)$ and those of $\sigma_\alpha$ on $Y$. In particular, if $\alpha(Y)$ intersects every orbit of $\sigma$ then $\sigma$ and $\sigma_\alpha$ have the same number of orbits.
\end{lem}
\begin{proof}
  Let $O$ be a $\sigma$-orbit of $X$, which by assumption contains $\alpha(y)\eqqcolon x$ for some $y\in Y$. It suffices to prove that $\alpha^{-1}(O)$ is a $\sigma_\alpha$-orbit. Now $O=\{\sigma^i(x):i\ge0\}$; let $0=m_0<m_1<\cdots$ be all the indices $m_i$ such that $\sigma^{m_i}(x)\in\alpha(Y)$, so $O\cap\alpha(Y)=\{\sigma^{m_i}(x):i\ge0\}$, and $\sigma^{m_i}(x)=\alpha(\sigma_\alpha^i(y))$, so $\alpha^{-1}(O)=\{\sigma_\alpha^i(y):i\ge0\}$.
\end{proof}

Consider now $1\le i<j<k<\ell\le n$, and denote by $\alpha$ the injective map $\{1,\dots,n-4\}\to\{1,\dots n\}$ defined by
\[1,\dots,n-4 \mapsto 1,\dots,i-1,k+1,\dots,\ell-1,j+1,\dots,k-1,i+1,\dots,j-1,\ell+1,\dots,n.\]
\begin{prop}\label{prop:ijkl}
  Assume $n\ge5$ and let $\pi\in S_n$ map $i,j,k,\ell$ to $k,\ell,i,j$, namely its cycle decomposition contains the cycles $(i,k)$ and $(j,\ell)$. Then $v(\pi)=v(\pi_\alpha)$.
\end{prop}
\begin{proof}
  Write $\rho=(1,\dots,n)(i,k)(j,\ell)$ and $\pi''=(i,k)(j,\ell)\pi$, the restriction of $\pi$ to the range of $\alpha$, so $\rho\pi''=(1,\dots,n)\pi$. An easy calculation checks $\rho_\alpha=(1,\dots,n-4)$. We note that the image of $\alpha$ intersects all orbits of $(1,\dots,n)\pi$: indeed it suffices to show that no orbit is contained in $\{i,j,k,\ell\}$. Now the images of $i,j,k,\ell$ are respectively $k+1,\ell+1,i+1,j+1\pmod n$ so this may happen only if $i+1=j,j+1=k,k+1=\ell,\ell+1=i-n$ and thus $n=4$. We thus have
  \begin{align*}
    v(\pi) &= \text{number of orbits of }(1,\dots,n)\pi\\
           &= \text{number of orbits of }((1,\dots,n)\pi)_\alpha\text{ by Lemma~\ref{lem:if2}}\\
           &= \text{number of orbits of }(\rho\pi'')_\alpha=\text{number of orbits of }\rho_\alpha\pi_\alpha\text{ by Lemma~\ref{lem:if1}}\\
           &= \text{number of orbits of }(1,\dots,n-4)\pi_\alpha = v(\pi_\alpha).\qedhere
  \end{align*}
\end{proof}

\subsection{Proof of Theorem~\ref{thm:main}}
We recall our main result:
\begin{thm}
  Consider a non-trivial word $w\in [F,F]$. Then there exists a factorization of $w$ without cancellations
  \[w=w_1a^{-1}w_2b^{-1}w_3a w_4b w_5\text{ with $a,b$ generators or their inverses}\]
  and $\|w_1w_4w_3w_2w_5\|=\|w\|-1$. Furthermore we have
  \begin{equation}\label{eq:wrecursion}
    w=[(w_4 w_3 a)^{w_1^{-1}},(b w_2^{-1} w_3^{-1})^{w_4^{-1}w_1^{-1}}]\cdot(w_1w_4w_3w_2w_5).
  \end{equation}
\end{thm}
\begin{proof}
  The last equality is directly checked by expanding the commutator; therefore $\|w_1w_4w_3w_2w_5\|\ge\|w\|-1$, and it remains to prove the reverse inequality. Write $n=|w|$, and let $\pi\in S_n$ be a pairing on $w$ maximizing $v(\pi)$. If $n=4$ then $w_1=w_2=w_3=w_4=w_5=1$ and the result follows, so we may assume $n\ge5$. Say two transpositions $(i,k),(j,\ell)$ with $i<k,j<\ell$ are \emph{linked} if the intervals $[i,k]$ and $[j,\ell]$ are neither nested nor disjoint. If all transpositions in the cycle decomposition of $\pi$ are unlinked, then $w$ freely reduces to the trivial word; we may thus assume that there are linked transpositions $(i,k),(j,\ell)$ in $\pi$, ordered so that $i<j<k<\ell$. By Proposition~\ref{prop:ijkl} and twice Theorem~\ref{thm:bardakov} we have
  \[\|w_1w_4w_3w_2w_5\|\le\frac{2-2v(\pi_\alpha)+n-4}4=\frac{2-2v(\pi)+n-4}4=\frac{2-2v(\pi)+n}4-1=\|w\|-1.\qedhere\]
\end{proof}

\begin{proof}[Proof of Corollary~\ref{cor:logspace}]
  We unroll $g$ times the recursion in Algorithm~\ref{algo:1}, arriving at an iterative algorithm with $4g$ nested loops.
  We implement it using a Turing machine whose read-write storage contains $4g$ pointers in the range $\{1,\dots,n\}$. Using these $4g$ pointers $i_1,j_1,k_1,\ell_1,\dots,i_g,j_g,k_g,\ell_g$, it is possible with $\mathcal O(\log n)$ memory to compute the function mapping $i\in\{1,\dots,n\}$ to the $i$th letter of the word obtained at the last stage of the recursion. We then apply the Lipton-Zalcstein algorithm to check in \textsc{LogSpace} whether this word is trivial.
\end{proof}

\begin{algorithm}[H]
  \KwData{A word $w$ and an integer $g\ge0$}
  \KwResult{A list of pairs of words expressing $w$ as a product of $g$ commutators, or \textsf{fail} if no such factorization exists}
  cyclically reduce $w$\;
  \If{$g=0$}{\Return \textsf{fail} \textbf{\textup{if}} $w\ne1$, \textbf{\textup{else}} $()$}
  \For{$1\le i<j<k<\ell\le|w|$}{
    \If{$w[i]\ne w[k]^{-1}$ \textbf{\textup{or}} $w[j]\ne w[\ell]^{-1}$}{\textbf{continue}}
    $(w_1,w_2,w_3,w_4,w_5)\leftarrow(w[1:i-1],w[i+1:j-1],w[j+1:k-1],w[k+1:\ell-1],w[\ell+1:|w|])$\;
    $F\leftarrow \operatorname{factorization}(w_1 w_4 w_3 w_2 w_5,g-1)$\;
    \If{$F\ne\textsf{\textup{fail}}$}{\Return $((w_1 w_4 w_3 w_1^{-1},w_1 w_4 w_2 w_3 w_4^{-1}w_1^{-1}),F)$}
    \Return \textsf{fail}
  }
  \caption{$\operatorname{factorization}(w,g)$}\label{algo:2}
\end{algorithm}

\begin{proof}[Proof of Corollary~\ref{cor:algo}]
  To compute the commutator length of $w$, we apply Algorithm~\ref{algo:1} first with $g=0$, then $g=1$, then $g=2$ etc. till it succeeds. In the RAM model, it is possible to keep track of $w_1w_4w_3w_2w_5$ as a linked list of letters, and at each elementary step of the algorithm (increase $i,j,k$ or $\ell$) a bounded number of operations need be executed to adjust the product $w_1w_4w_3w_2w_5$, including its free cancellations. Therefore the complexity of the algorithm, at each step of the recursion, is controlled by the loop over $i<j<k<\ell$, with time complexity $\mathcal O(|w|^4)$.

  Furthermore, Algorithm~\ref{algo:1} may be modified in a straightforward manner to return an expression of $w$ as a product of commutators, see Algorithm~\ref{algo:2}.
\end{proof}

\subsection{Implementation}\label{ss:implementation}
It is possible to speed up a little bit the algorithm, and this was crucial for the computations in the next section. Let us concentrate on the case $\|w\|=1$, namely determine whether $w$ is itself a commutator. A naive implementation, looping over all positions of $a^{-1},b^{-1},a,b$ and testing whether $w_1w_4w_3w_2w_5$ freely reduces to the trivial word, requires $\mathcal O(|w|^5)$ steps: four factors $|w|$ for the loop and one for the reduction of $w_1w_4w_3w_2w_5$.

This time can be reduced to $\mathcal O(|w|^3)$ assuming large machine integers, and even $\mathcal O(|w|^2)$ on average, as follows: firstly, rather than keeping track of words $w_1,w_2,w_3,w_4,w_5$ we store $2\times2$ integer matrices faithfully representing them. If $F$ has rank $2$, we can even choose as images of generators the transvections $(\begin{smallmatrix}1&2\\0&1\end{smallmatrix})$ and $(\begin{smallmatrix}1&0\\2&1\end{smallmatrix})$. We update the words $w_i$ at each step of the algorithm by elementary row and column operations. In this manner, words of length up to $64$ may faithfully be stored into four $64$-bit integers. Next, we note that we may loop first on the positions of $a^{-1}$ and $a$; and that, if $w$ is a commutator, then the words $w_3^{-1} b w_2^{-1}$ and $w_4 b w_5w_1$ are conjugates of each other. Knowing their matrix representations, we can immediately rule out a pair of $a$-positions if these matrices have different traces, a constant-time test that eliminates almost all candidates (note that integer overflow is not an issue here, since we are looking for cheap ways of ruling out some candidates).

For those that survive the test, rather than considering the $\mathcal O(|w|^2)$ pairs of $b$-positions, we could cyclically reduce these words and check whether they are cyclic permutations of each other, using for example the Knuth-Morris-Pratt algorithm~\cite{kmp:matching}; but this case occurs so seldom in practice that it was not worth implementing.

In even more detail: our implementation loops over $i,k,j,\ell$ in that order, and maintains matrices $M_{23451},M_{23},M_{451},M_{32},M_{514}$ representing the respective products $w_2 b^{-1}w_3 a w_4 b w_5w_1,$ $w_2 b^{-1}w_3,$ $w_4 b w_5 w_1,$ $w_3w_2,$ $w_5w_1w_4$. At each elementary step (increase one of $i,j,k,\ell$), one of these matrices is to be multiplied on the left or on the right by a generator. For each choice of $i,k$, the implementation checks whether $M_{451}$ and $M_{23}$ have same trace (recall that in $\operatorname{SL}_2(\mathbb Z)$ the trace of a matrix equals the trace of its inverse), before looping on $j$ and $\ell$.

\section{Quadratic equations and Theorem~\ref{thm:32}}\label{ss:quadratic}
Consider the equation $x_1^2x_2^2x_3^2x_4^2=1$ in a free group $F$. A solution $(x_1,x_2,x_3,x_4)$ is the same thing as a homomorphism $\phi\colon S\to F$, with $S=\langle x_1,x_2,x_3,x_4\mid x_1^2x_2^2x_3^2x_4^2=1\rangle$. According to~\cite{piollet:quadratique}, the set $\Phi$ of solutions $\phi$ to the equation is characterized as follows: let $F_2$ denote the free group of rank $2$. A \emph{literal} solution $\lambda\colon S\to F_2$ is a homomorphism sending each generator to a word of length $\le1$. The \emph{mapping class group} of $S$ is the group of automorphisms of $S$ that preserve $\{x_1\}^S\cup\cdots\cup\{x_4\}^S$; equivalently, automorphisms of the free group $\langle x_1,\dots,x_4\rangle$ that preserve the conjugacy class of $x_1^2x_2^2x_3^2x_4^2$. Then
\begin{equation}\label{eq:Phi}
  \Phi=\{\eta\circ\lambda\circ\alpha:\alpha\in\operatorname{MCG}(S),\,\lambda\text{ literal},\,\eta\colon F_2\to F\text{ any homomorphism}\}.
\end{equation}

\subsection{The mapping class group of the non-orientable surface of genus $2$}
In order to explore the solutions of the equation $x_1^2x_2^2x_3^2x_4^2=1$, it is helpful to also consider a different presentation of $S$, in which generators of $\operatorname{MCG}(S)$ are easier to write. We use
\[S=\langle s_1,t_1,s_2,t_2\mid s_1t_1s_1^{-1}t_1^{-1}s_2t_2s_2^{-1}t_2\rangle.\]
The isomorphism is given for instance by
\[(s_1,t_1,s_2,t_2)=(x_2x_3^2x_4,x_3^{-1}x_2^{-1},x_3x_4x_2,x_3x_4x_1x_2).\]
Using these generators, we have the following mapping classes, in which unlisted generators are fixed:
\begin{gather*}
  \alpha_1\colon\begin{cases} t_1 &\mapsto t_1s_1^{-1}\end{cases}\qquad
  \beta_1\colon\begin{cases} s_1 &\mapsto s_1t_1,\end{cases}\qquad
  \beta_2\colon\begin{cases} s_2 &\mapsto s_2t_2,\end{cases}\\
  \gamma_1\colon\begin{cases} s_1 &\mapsto s_1u\\ t_1 &\mapsto u^{-1}t_1u\\ s_2 &\mapsto u^{-1}s_2\end{cases}\text{ for }u=s_2t_2^{-1}s_2^{-1}t_1,\\
  \eta_3\colon\begin{cases}
    s_1 &\mapsto s_1v\\
    t_1 &\mapsto v^{-1}t_1v \text{ for }v=s_2^2t_2^{-1}s_2t_2^{-1}s_2^{-1}t_1,\\
    s_2 &\mapsto v^{-1}t_1s_2\\
    t_2 &\mapsto t_1^{-1}s_2t_2s_2^{-2}t_1s_2,\\
  \end{cases}
\end{gather*}
and by~\cite{comerford-lee:2comm}*{Appendix} these elements generate $\operatorname{MCG}(S)$.

We add for good measure the automorphism $\delta\colon x_1\mapsto x_2\mapsto x_3\mapsto x_4\mapsto x_1$ in the original generators, and thus consider as generating set for $\operatorname{MCG}(S)$ the collection
\begin{equation}\label{eq:Sigma}
  \Sigma=\{\alpha_1^{\pm1},\beta_1^{\pm1},\beta_2^{\pm1},\gamma_1^{\pm1},\eta_3^{\pm1},\delta,\delta^2,\delta^3\}.
\end{equation}
In the generating set $\{s_1,t_1,s_2,t_2\}$, there are three maximal literal solutions $\tau_{-1},\tau_0,\tau_1$ given for $i=-1,0,1$ by
\[\tau_i\colon\begin{cases}s_1 &\mapsto x,\\ t_1 &\mapsto x^i,\\ s_2 &\mapsto y,\\ t_2 &\mapsto 1.\end{cases}\]
  
\subsection{Proof of Theorem~\ref{thm:32}}
The connection to Theorem~\ref{thm:32} is the following: given a solution $(x_1,\dots,x_4)\in F^4$ to $x_1^2x_2^2x_3^2x_4^2=1$, consider the element $w=x_1x_2x_3x_4$. Then
\begin{align*}
  w^2 &= (x_1x_2x_3x_4)^2 (x_1^2x_2^2x_3^2x_4^2)^{-1}\\
  &= \big([x_4^{-1}x_3^{-1},x_3^{-1}x_2^{-1}x_1^{-1}]\cdot[x_2^{-1}x_1^{-1},x_2]\big)^{x_2^{-1}x_1^{-1}},
\end{align*}
so for every solution $(x_1,x_2,x_3,x_4)$ the element $w^2=(x_1x_2x_3x_4)^2$ has commutator length at most $2$. On the other hand, the first claim of Theorem~\ref{thm:32} will follow by constructing a sufficiently complicated solution.

We considered, in~\eqref{eq:Phi}, all $\alpha$ a product of at most $5$ of the generators from $\Sigma$ given in~\eqref{eq:Sigma}; all $\lambda\in\{\tau_{-1},\tau_0,\tau_1\}$; and $F=F_2$ and $\eta=1$. We computed, for each solution, the image $w$ of $x_1x_2x_3x_4$ (which is conjugate to $t_2$), and sorted them in increasing order of length. We finally computed the commutator length of each of these solutions (see the next section), and after examining about 2500 candidates arrived at a solution
\begin{multline}\label{eq:w}
  w=x^{-2} y x y^2 x^{-2} y x y x^{-1} y^{-1} x y^2 x^{-1} y^2 x y^{-1} x^{-1} y^{-1} x^2 y^{-2} x^{-1} y^{-1} x y^{-2}\\\cdot x y^{-2} x^{-1} y^{-1} x y^{-2} x^{-1} y x y^{-1} x^{-1} y^{-1} x^2 y^{-1} x^{-1} y^{-1} x y^2 x^{-1} y x y^2 x^{-2} y x y
\end{multline}
with commutator length $3$. It was produced by $\lambda=\tau_1$ and $\alpha=\eta_3^{-1}\circ\alpha_1\circ\eta_3^{-1}\circ\delta^2\circ\eta_3^{-1}$. Dozens of other solutions appeared, but no shorter one.

Algorithm~\ref{algo:1}, with the improvements described in~\S\ref{ss:implementation}, could test words of length comparable to $w$ in about a second. Words of length above 100 could be routinely tested. For extra safety and to protect against hidden bugs, the commutator length of $w$ was also computed using Alden Walker's program \texttt{scallop}~\cite{walker:scallop}. With the options `\texttt{-cyclic -C}' and the \texttt{Gurobi} solver~\cite{gurobi:web} as backend, it certified $\|w\|=3$ for the word $w$ from~\eqref{eq:w} in a bit more than one hour.

\subsection{Open problems}
There does not seem to be much hope to obtain, by brute force, words $w$ with $\|w\|\ge4$ and $\|w^2\|=2$. However, the construction of $w$ in~\eqref{eq:w} in the form $\tau_1(\alpha(x_1x_2x_3x_4))$ for a mapping class $\alpha$ that is a product of $5$ generators raises the following problem:
\begin{question}
  Does there exist a mapping class $\alpha\in\operatorname{MCG}(S)$ with $\|\tau_1(\alpha^n(x_1x_2x_3x_4))\|\ge n$ for all $n\in\N$?
\end{question}
This would provide a systematic collection of solutions, and remove much of the nonconstructivity of the implicit function theorem in free groups.

It is as yet unknown whether commutator length can decrease upon taking \emph{cubes}; more importantly, whether there exist a sequence $(w_n)$ in $F$ with $\|w_n\|\to\infty$ and $\|w_n^3\|$ bounded, or equivalently whether there is $3$-torsion in the abelianization of $F^\infty$. Following the same idea, one would want arbitrarily complicated solutions to the equation $x_1^3x_2^3x_3^3x_4^3$, a goal that seems out of reach now.

Danny Calegari suggested a variant of the question, that may be more tractable: ``can one find $w$ with $\|3w\|<\|w\|$?''. Here $\|m\cdot w\|$ is the minimal number of commutators required to express a product of $m$ conjugates of $w$, and is trivially at most $\|w^m\|$. A sequence $(w_n)$ with $\|w_n\|\to\infty$ and $\|3w_n\|$ bounded would likewise imply the existence of $3$-torsion in $H_1(F^\infty)$.


\end{document}